\documentclass[12pt]{article}

\usepackage[ruled]{algorithm2e}
\usepackage{pdfsync}
\usepackage{graphicx}
\usepackage{subfigure}
\usepackage{epsfig}
\usepackage{booktabs}
\usepackage{amsmath}
\usepackage{amsfonts}
\usepackage{cite}
\usepackage{url}
\usepackage[margin=1in]{geometry}

\DeclareGraphicsExtensions{.pdf}

\title{MATSuMoTo: The MATLAB Surrogate Model Toolbox For Computationally Expensive   Black-Box Global Optimization Problems}
\author{Juliane M\"uller\thanks{Corresponding author: Cornell University, School of Civil and Environmental Engineering, 220 Hollister Hall, Ithaca, NY 14853-3501; {\href{mailto:juliane.mueller2901@gmail.com}{\tt juliane.mueller2901@gmail.com}}}}
\providecommand{\keywords}[1]{\textbf{\textit{Index terms---}} #1}

\begin{document}

\maketitle

\begin{abstract}
MATSuMoTo is the MATLAB Surrogate Model Toolbox  for computationally expensive, black-box, global optimization problems that may have continuous, mixed-integer, or pure integer variables. Due to the black-box nature of the objective function, derivatives are not available. Hence, surrogate models are used as computationally cheap approximations of the expensive objective function in order to guide the search for improved solutions. Due to the computational expense of doing a single function evaluation, the goal is to find optimal solutions within very few expensive evaluations.  The multimodality of the expensive black-box function requires an algorithm that is able to search locally as well as globally. MATSuMoTo is able to address these challenges. MATSuMoTo offers  various choices for surrogate models and surrogate model mixtures,  initial experimental design strategies, and  sampling strategies. MATSuMoTo is able to do several function evaluations in  parallel by exploiting MATLAB's Parallel Computing Toolbox.     
\end{abstract}

%
%
%
%
%
%

\maketitle

\section{Introduction and Motivation}
Simulation models  for approximating physical phenomena such as, for example,  cleaning up contaminated groundwater  or simulating the global climate are becoming with increasing computational power more complex.  Many  simulation models  require a substantial amount of computation time. One simulation may take from several minutes to several hours or even days. If parameters of these models have to be optimized with the goal, for example,   to minimize the error between simulation model predictions and observations,  the computational expense arising  from a single simulation restricts the total number of  function evaluations that can be allowed in order to obtain a  solution within reasonable time. Hence, algorithms that efficiently find a near optimal solution are needed.\\

Furthermore, these simulation models are often black boxes and an analytical description of the objective function is not available. Therefore, derivative information is not available. Automatic differentiation can in many cases not be applied due to the large number of operations required in the simulation or due to confidentiality restrictions of the simulation code. Numerically approximating derivatives, on the other hand, requires too many expensive function evaluations and is therefore not feasible with respect to computation times.\\

Since the objective function is a black box, it is not known whether or not multiple local and global optima exist. If it is known that the problem  is one of local optimization (in this case the local optimum is also the global optimum), a local optimization algorithm should be applied. However,  in the case of multimodal optimization problems, local optimization algorithms stop generally at a local optimum and are not able to find the global optimum. Thus, for black-box optimization problems for which it is not known how many local and global optima exist, it is necessary to use a global optimization algorithm that is able to continue the search globally after a local optimum has been detected.  \\


In this paper we address  optimization  problems of the following type:
\begin{equation}
\min_{\mathbf{x}\in\mathcal{D}} f(\mathbf{x}),
\end{equation} 

where $\mathbf{x}$ is the variable vector,  $f(\mathbf{x})$ represents a continuous, deterministic, computationally expensive, black-box objective function, and $\mathcal{D}$ denotes the  box-constrained variable domain defined by $-\infty <x_i^{\text{low}} \leq x_i \leq x_i^{\text{up}} < \infty$, $i = 1,\ldots,d$, where $d$ denotes the problem dimension. The variable vector $\mathbf{x}$ may consist of 
\begin{itemize}
\item only continuous variables, i.e.\ $\mathbf{x} \in \mathbb{R}^{d_1}$, where $d_1$ is the  dimension of the continuous variables,
\item  continuous and integer variables, i.e.\ $\mathbf{x} = (\mathbf{u}^T, \mathbf{v}^T)^T$ where $\mathbf{v} \in \mathbb{R}^{d_1}$ and $\mathbf{u} \in\mathbb{Z}^{d_2}$, where $d_1$ and $d_2$ are the number of continuous and integer variables, respectively,
\item only integer variables, i.e.\ $\mathbf{x} \in \mathbb{Z}^{d_2}$, where $d_2$ is the  dimension of the integer variables.
\end{itemize}

Derivative-free algorithms~\cite{Conn2009} have been developed for addressing problems for which derivative information is not available.  In this paper we present MATSuMoTo, the \underline{MAT}LAB \underline{Su}rrogate \underline{Mo}del \underline{To}olbox, which  uses, in particular, surrogate models~\cite{Booker1999} for tackling the challenge of finding (near) optimal solutions of computationally expensive black-box problems  within a very limited number of function evaluations. The surrogate model approach requires in comparison to population-based algorithms such as, for example, genetic algorithms and particle swarm methods, significantly fewer function evaluations and is therefore better suitable for problems where the function evaluation is computationally very expensive~\cite{Muller2013, Muller2012, Muller2014, Regis2007b, Wild2007}.\\

We developed MATSuMoTo for optimization problems where the function evaluation requires a significant amount of computation time (several minutes or more). For such problems, the computational overhead from building the surrogate model is insignificant compared to evaluating the objective function. For problems where the objective function can be cheaply evaluated (fractions of a second), we do not recommend using MATSuMoTo because in that case building the surrogate model will be the driver for the computation time. \\

MATSuMoTo requires MATLAB version 2010b or newer. The MATLAB Statistics Toolbox is required. If the user wishes to make use of MATSuMoTo's option for simultaneously evaluating the objective function at various points, the MATLAB Parallel Computing Toolbox is required. Depending on the sampling strategy,  MATLAB's Global Optimization Toolbox and Optimization Toolbox are required. The algorithm can be used without the Parallel Computing and Optimization toolboxes, but in that case not all features will be available (see the user manual for the dependencies of MATSuMoTo's features on MATLAB toolboxes).  A Python implementation of the toolbox is in development. To our best knowledge, MATSuMoTo  is the only toolbox that is able to address the variety of optimization problems described above  using surrogate model algorithms with the rich selection of features described in this paper. \\

MATSuMoTo is available on  GitHub\footnote{\url{https://github.com/Piiloblondie/MATSuMoTo}} and the author's website\footnote{\url{https://courses.cit.cornell.edu/jmueller/}}. The toolbox comes with a thorough  user manual that describes how to use MATSuMoTo and explains  the various options  in detail. MATSuMoTo contains features of the algorithms described in~\cite{Muller2010, Muller2012, Muller2013, Muller2014}.\\

This paper is organized as follows.  In Section~\ref{sec:surrogates}, we briefly summarize what surrogate models are and what types of surrogate models exist. The steps of a general surrogate model algorithm are described and illustrated in Section~\ref{sec:surrogatealgorithms}. The implementation of MATSuMoTo and the various options for the elements of the surrogate model algorithm (initial experimental design, type of surrogate model, sampling strategy) are detailed  in Section~\ref{sec:MATSuMoTo}. Section~\ref{sec:discussion} summarizes the paper and outlines ongoing  developments of features to be included in MATSuMoTo in the future.

\section{Surrogate Models} \label{sec:surrogates}
Computationally expensive simulation models are used to approximate complex physical phenomena. Similarly, surrogate models are used to approximate the computationally expensive simulation model. Surrogate models (also known as response surface models or metamodels) are used in place of  the expensive simulation model during the optimization~\cite{Booker1999}, i.e.\ $f(\mathbf{x}) = s(\mathbf{x}) + e(\mathbf{x})$, where $f(\mathbf{x})$ represents the output of the computationally expensive simulation model, $s(\mathbf{x})$ is the prediction of the surrogate model at point $\mathbf{x}$, and $e(\mathbf{x})$ denotes the difference between both. During the optimization, information from the computationally cheap surrogate model is used to guide the search for promising points in the variable domain, and thus the expensive objective function is evaluated only at carefully selected points. Hence, the total number of required function evaluations to find (near) optimal solutions is significantly decreased and the total optimization time is reduced.  \\

Various surrogate model types  have been developed in the literature and applied to a large variety of application problems. There are interpolating models, such as kriging~\cite{Jones1998, Jones2001, Lam2009a, Martin2005, Matheron1963}	 and radial basis function models~\cite{Gutmann2001, Powell1992}, and there are non-interpolating models such as multivariate adaptive regression splines~\cite{Friedman1991} and polynomial regression models. Polynomial regression models have widely been used in engineering applications~\cite{Myers1995}. Recently, radial basis function models have often been used for solving computationally expensive optimization applications~\cite{Bjork2001, Gutmann2001, Muller2013, Muller2012, Regis2007b, Regis2009, Regis2013, Wild2013}. Moreover, ensembles of surrogate models, i.e.\ surrogate models whose prediction is a weighted sum of predictions of individual models, have been developed~\cite{Goel2007, Muller2010, Viana2008}. The question of which surrogate model performs best and should be chosen is addressed in\cite{Muller2014}.

\section{Surrogate Model Algorithms}\label{sec:surrogatealgorithms}
\begin{figure}[h]
\centering
\includegraphics[width=10cm, height=10cm]{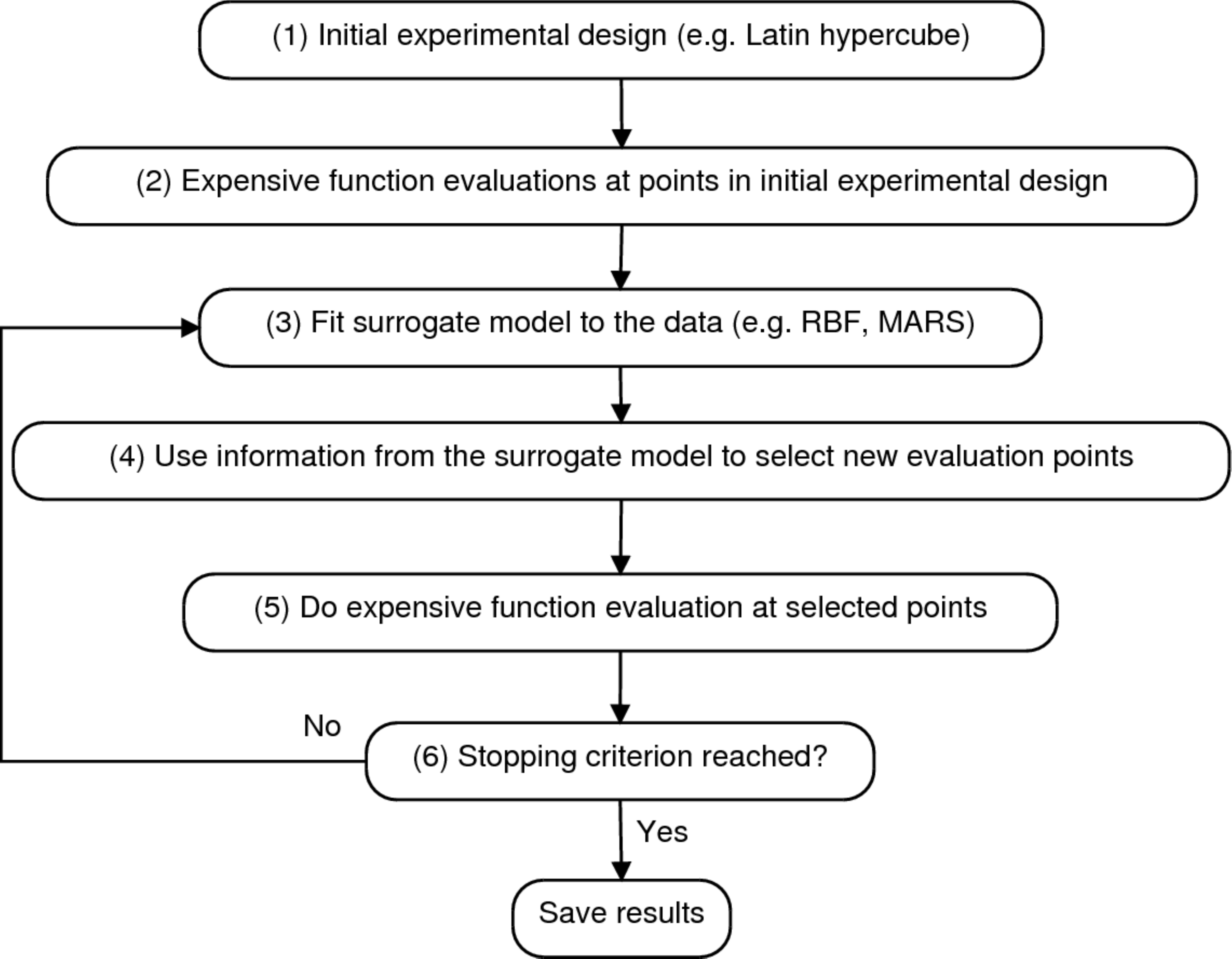}
\caption{Steps of a general surrogate model algorithm. In every iteration one or more points are selected for doing the expensive objective function evaluation and the surrogate model is updated with the new data.}
\label{Fig:AlgorithmSteps}
\end{figure}

Surrogate model algorithms consist in general of the steps illustrated in Figure~\ref{Fig:AlgorithmSteps}.
%
For generating the initial experimental design in Step~(1), any strategy may be used, for example, a (symmetric) Latin hypercube design, orthogonal arrays, or using (some or all of) the corner points of the variable domain (see~\cite{Montgomery2001} for various design strategies). It must be ensured that the number of selected points  suffices to fit the chosen surrogate model. In Step~(2), the computationally expensive black-box simulations have to be done at the selected design points. \\

In Step~(3), the surrogate model is fit to  the data. When using radial basis functions (RBF) and polynomial regression models, model parameters are calculated. Multivariate adaptive regression splines (MARS) are parameter-free and use a forward and backward pass to find a model that consists of a combination of  hinge functions that best approximates the data (see~\cite{Friedman1991}). The surrogate model is in Step~(4) used to predict the objective function value at unsampled points in the variable domain. This is computationally inexpensive because the true objective function is not evaluated in this step. Based on the predictions, a criterion is used to determine the most promising point(s) for doing the next expensive function evaluation(s). Such criteria might be, for instance, to use the minimum point of the surrogate model~\cite{Regis2012} or  to use a weighted score of the distance to already evaluated points and the objective function value prediction~\cite{Muller2014, Regis2007b}. \\

\begin{figure}[h!]
        \centering       
        \subfigure[Initial experimental design (black crosses).]{
                \includegraphics[width=0.45\textwidth, height=4cm]{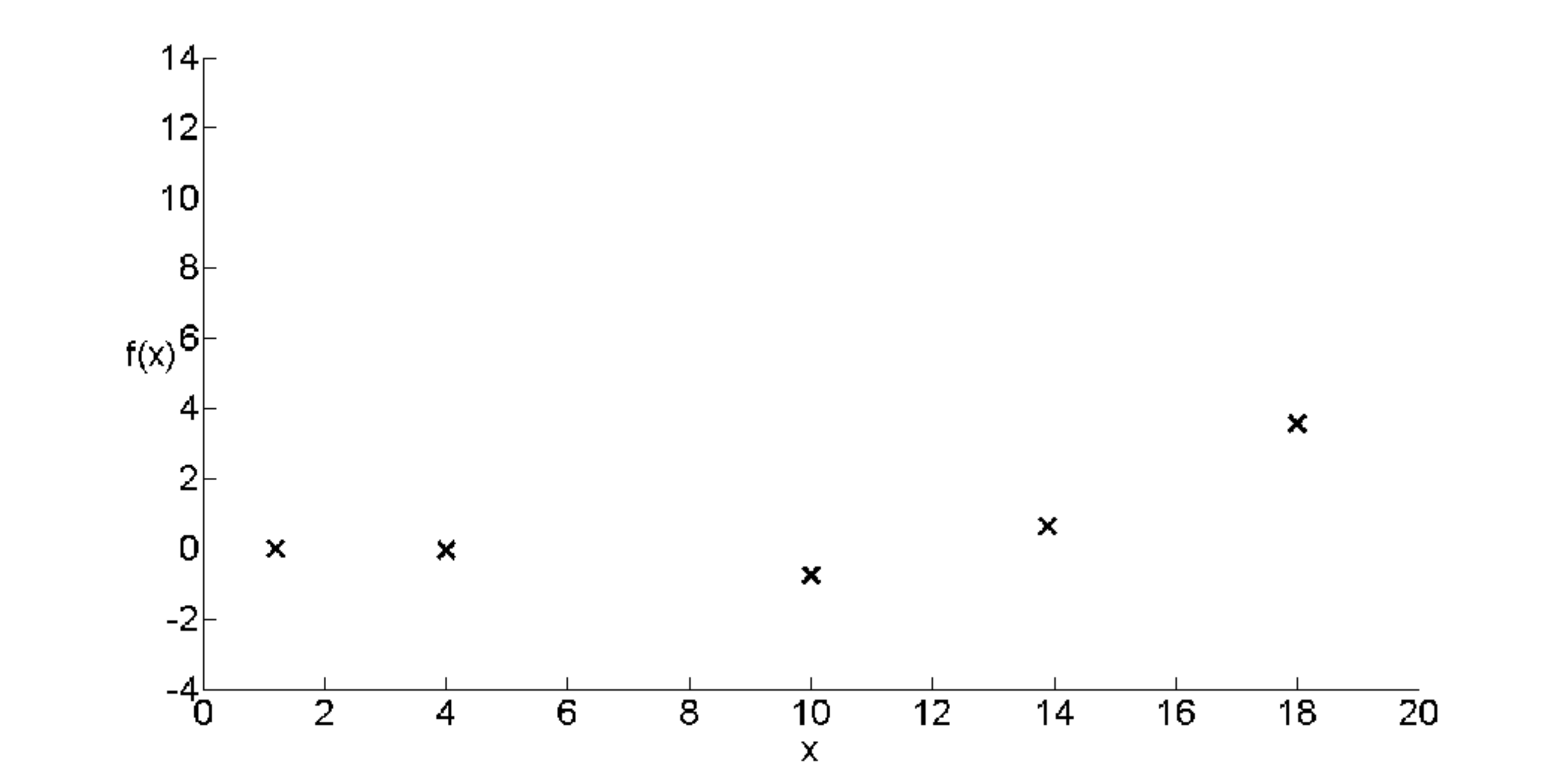}
                \label{SubFig: Tile_3}
        }
        \subfigure[Fit the  surrogate model (red graph).]{
                \includegraphics[width=0.45\textwidth, height=4cm]{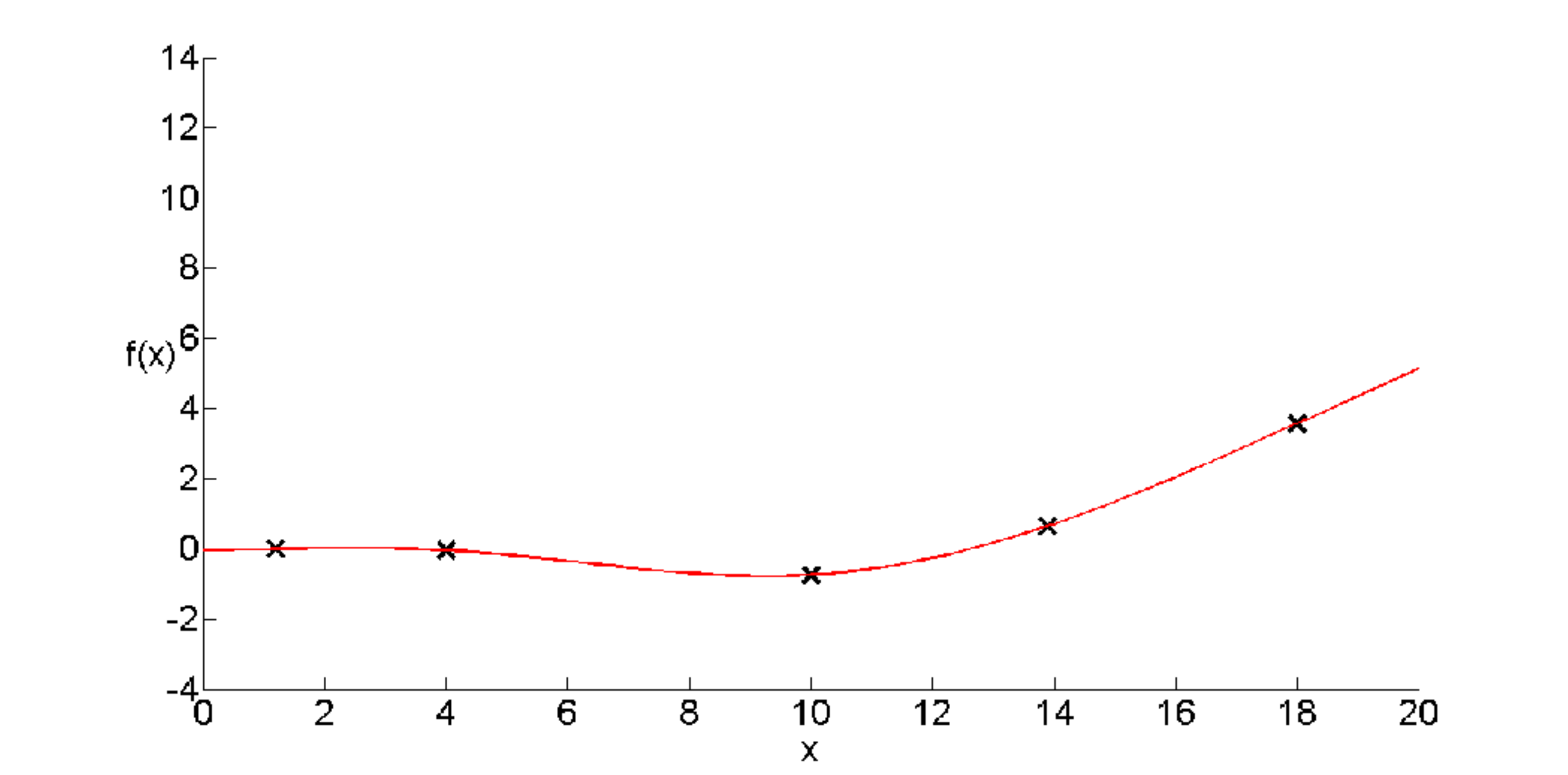}
                \label{SubFig: Tile_4}
        }
        
         \subfigure[Selected a new data point (blue circle).]{
                \includegraphics[width=0.45\textwidth, height=4cm]{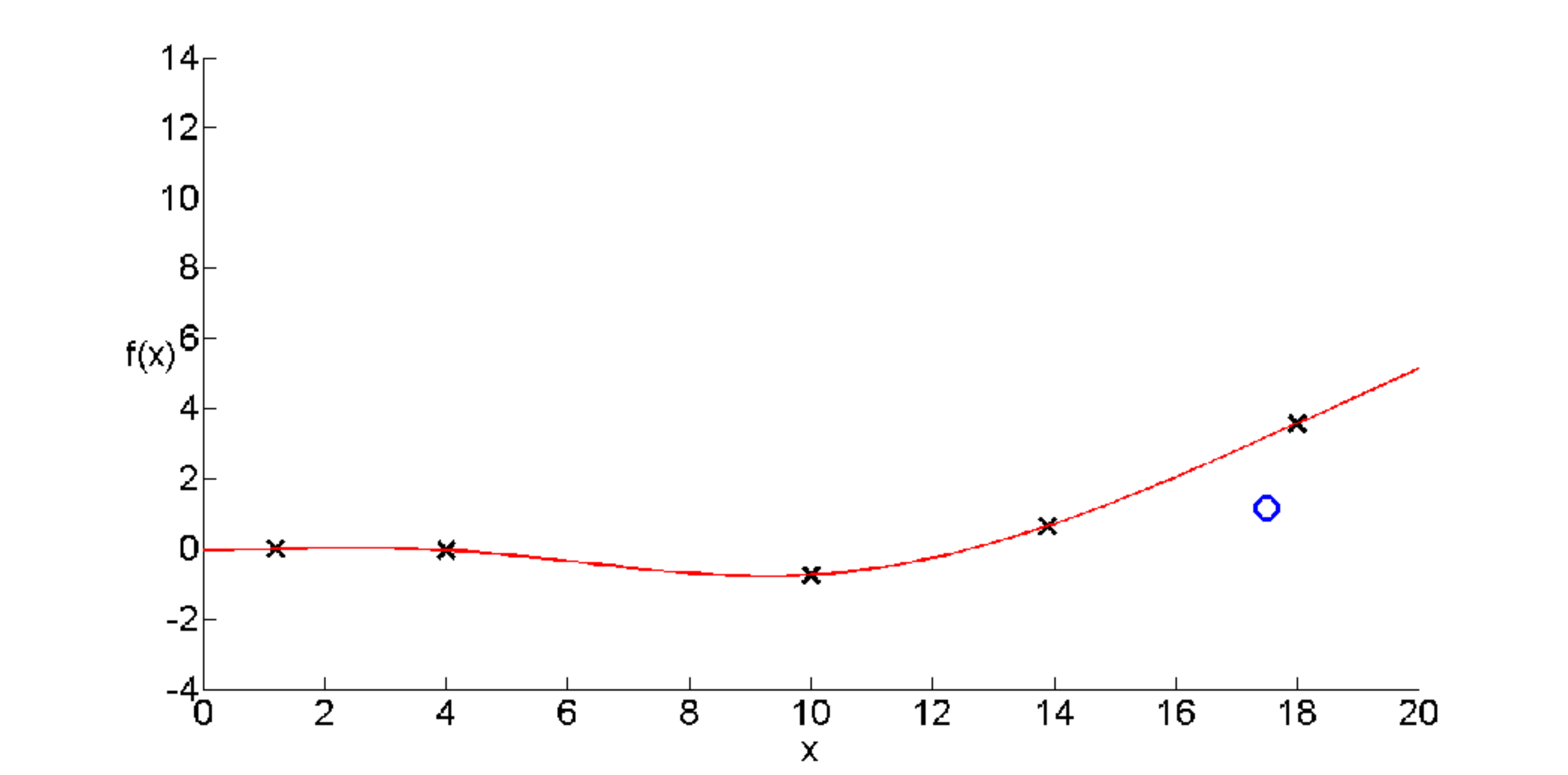}
                \label{SubFig: Tile_5}
        }
        \subfigure[Update the  surrogate model (red graph).]{
                \includegraphics[width=0.45\textwidth, height=4cm]{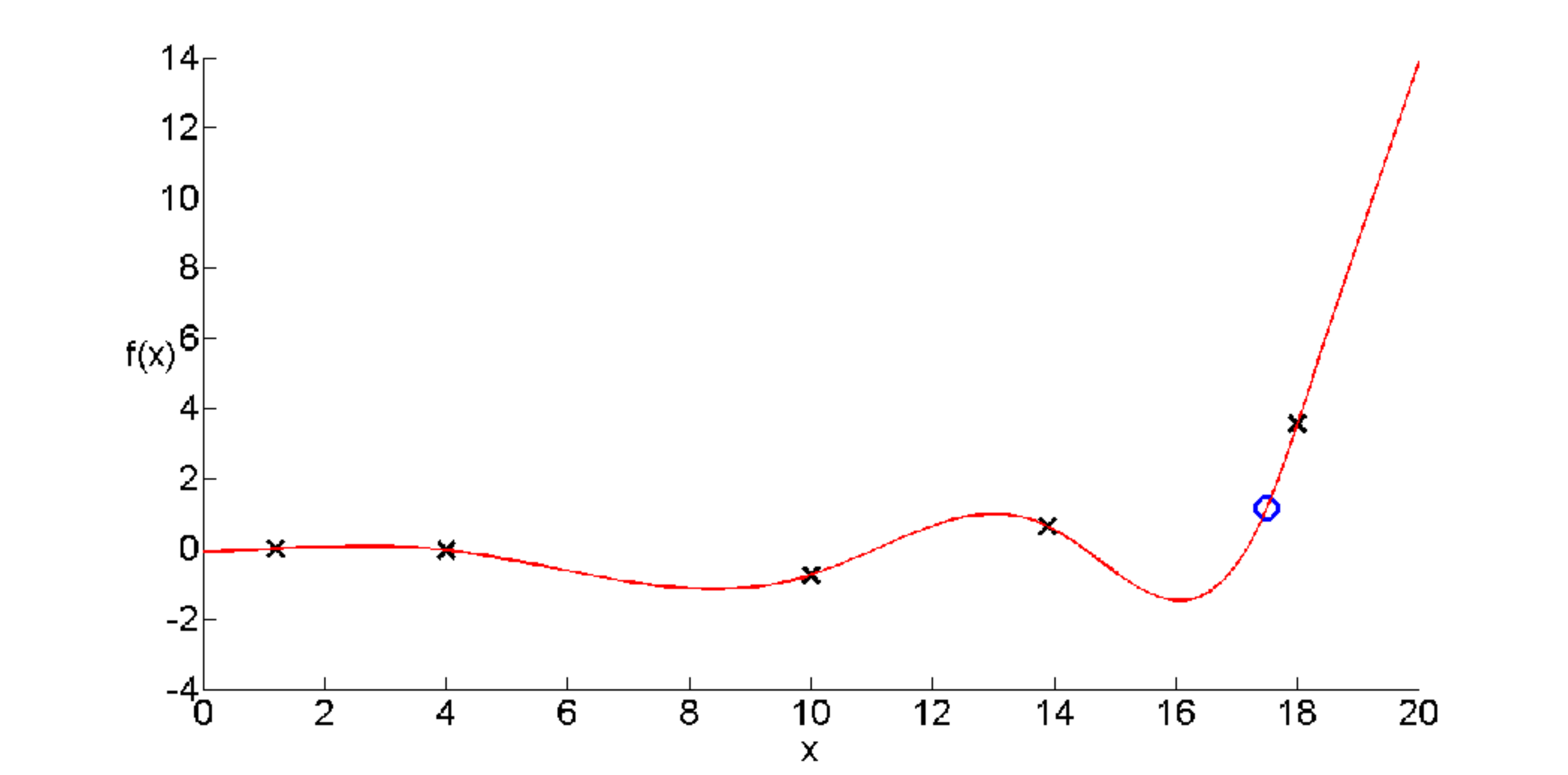}
                \label{SubFig: Tile_6}
        }
        
         \subfigure[Selected a new data point (blue circle).]{
                \includegraphics[width=0.45\textwidth, height=4cm]{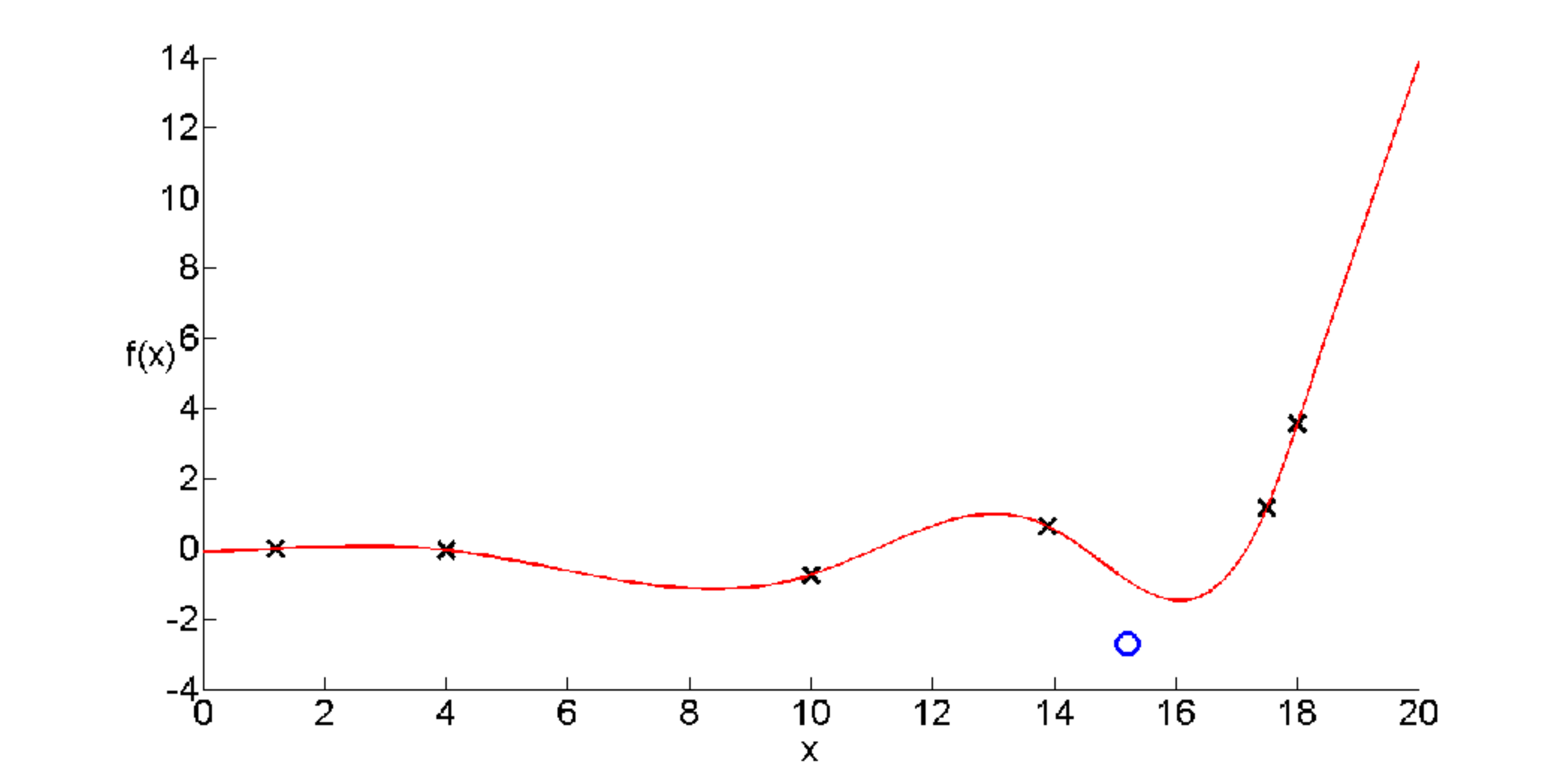}
                \label{SubFig: Tile_7}
        }
        \subfigure[Update the  surrogate model (red graph).]{
                \includegraphics[width=0.45\textwidth, height=4cm]{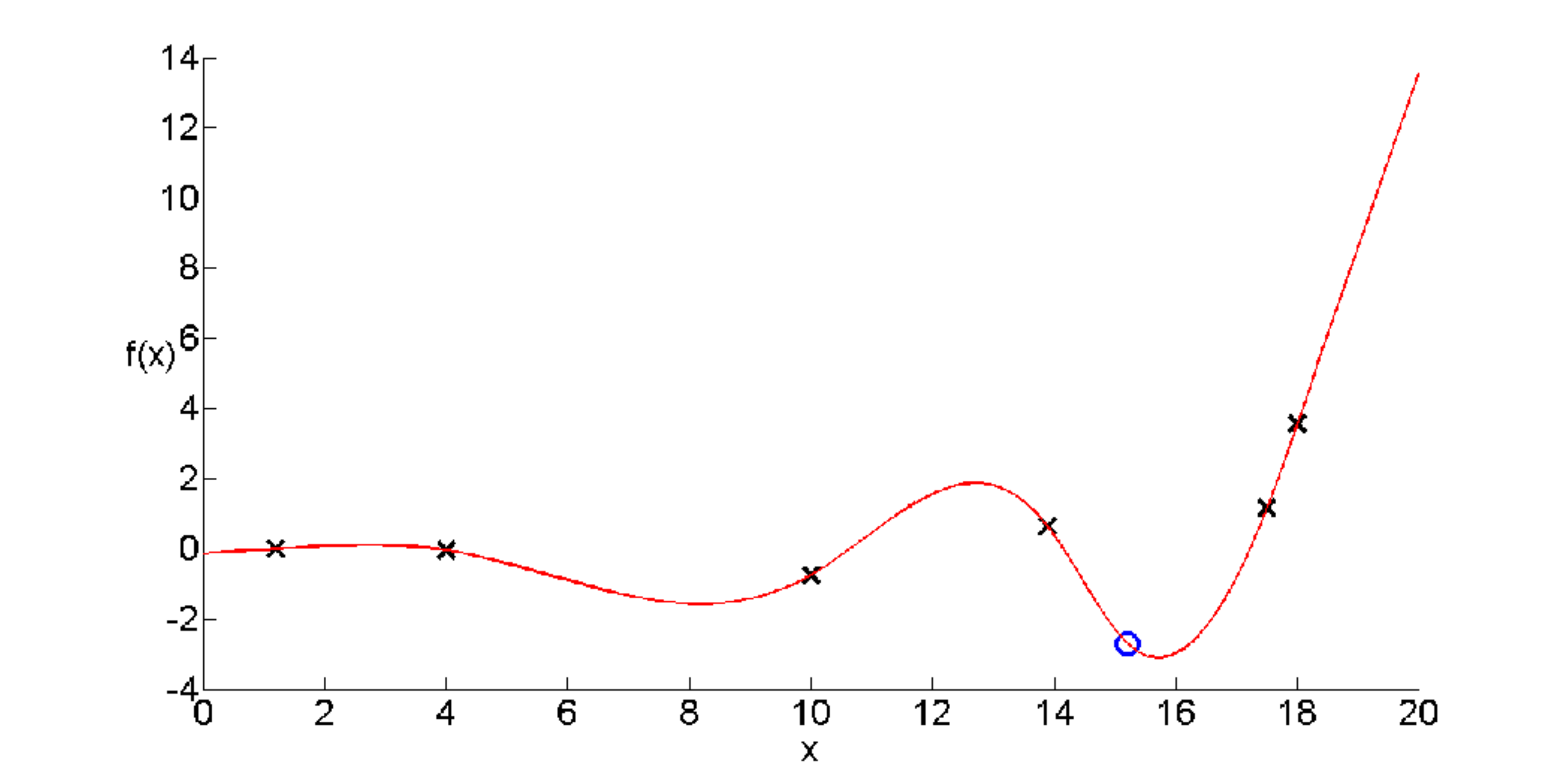}
                \label{SubFig: Tile_8}
        }
\caption{Illustration of the steps of the  surrogate model algorithm with a one-dimensional problem. Black crosses denote points that have already been evaluated in previous iterations, blue circles denote newly chosen data points, the red graph illustrates the interpolating  surrogate model.}
\label{Fig: SMalgo}
\end{figure}

The computationally expensive function  is evaluated at the selected point(s) in Step~(5) and  the new data is added to the set of already evaluated points. If the stopping criterion is not reached in Step~(6), the surrogate model is updated with the new data. The algorithm iterates between Steps~(3)-(6) until the stopping criterion has been met. Stopping criteria may be, for example, a maximal number of allowed function evaluations (used in MATSuMoTo), a maximal allowable CPU time, or a maximal number of failed consecutive improvement trials. Figure~\ref{Fig: SMalgo} illustrates the individual steps of the algorithm with a one-dimensional example. 

\section{MATSuMoTo Implementation} \label{sec:MATSuMoTo}
In this section we describe  the specific implementation of the individual steps (1)-(6) in MATSuMoTo and the most important toolbox features. Table~\ref{tab:settings} gives an overview over the choices that can be made. A more detailed description of the features can be found in the user manual.

\begin{table}[ht]
\centering
\caption{Possible feature choices for the individual steps of MATSuMoTo.}
\label{tab:settings}
\begin{tabular}{lll}
\toprule
Algorithm Step & Choice Name & Description\\
\midrule
(1) Initial design & \textsf{CORNER} & Corner points of the hypercube\\
& \textsf{SLHD} & Symmetric Latin hypercube\\
& \textsf{lhd} & Latin hypercube\\
\midrule
(3) Surrogate model & \textsf{RBFcub} & Cubic RBF\\
& \textsf{RBFtps} & Thin-plate spline RBF\\
& \textsf{RBFlin} & Linear RBF\\
& \textsf{MARS} & Multivariate adaptive regression spline\\
& \textsf{POLYlin}&Linear regression polynomial\\
& \textsf{POLYquad} & Quadratic regression polynomial\\
& \textsf{POLYquadr} & Reduced quadratic regression polynomial\\
& \textsf{POLYcub} & Cubic regression polynomial\\
& \textsf{POLYcubr} & Reduced cubic regression polynomial\\
& \textsf{MIX\_RcM}& Mixture of  \textsf{RBFcub} and \textsf{MARS} \\
& \textsf{MIX\_RcPc}& Mixture of  \textsf{RBFcub} and \textsf{POLYcub}\\
& \textsf{MIX\_RcPcr}& Mixture of \textsf{RBFcub} and \textsf{POLYcubr} \\
& \textsf{MIX\_RcPq}& Mixture of  \textsf{RBFcub} and \textsf{POLYquad}\\
& \textsf{MIX\_RcPqr}& Mixture of \textsf{RBFcub} and  \textsf{POLYquadr}\\
& \textsf{MIX\_RcPcM}&Mixture of  \textsf{RBFcub}, \textsf{POLYcub}, and \textsf{MARS}\\
\midrule
(4) Sampling strategy &  \textsf{CANDloc} & Local candidate point search\\
& \textsf{CANDglob} &Global candidate point search\\
& \textsf{SurfMin} &Minimum point of surrogate model\\
\bottomrule
\end{tabular}
\end{table}

The main file to start the algorithm is called \textsf{MATSuMoTo.m}. The user may give the following input arguments which are also described in more detail in the code manual:

\begin{enumerate}
\item[IA-1] the file describing the optimization problem (mandatory)
\item[IA-2] the maximum number of allowed expensive function evaluations (this is the algorithm's stopping criterion, optional)
\item[IA-3] the surrogate model type (optional)
\item[IA-4] the sampling strategy (optional)
\item[IA-5] the type of initial experimental design (optional)
\item[IA-6] the number of points  in the initial experimental design  (optional)
\item[IA-7] specific points to be included in the initial experimental design (optional)
\item[IA-8] the number of points to be selected in each iteration for doing the expensive function evaluation (optional)
\end{enumerate}

MATSuMoTo starts by validating the user's input and assigns default values to non-defined input arguments. MATSuMoTo then calls, depending on whether or not integer variables have been defined,  either the optimization routine for continuous, integer, or mixed-integer  problems. All three optimization routines work similarly and differ mainly with respect to the creation of the candidate points for doing the next expensive function evaluation (see Section~\ref{sec:sampling}).\\
 
\subsection{IA-1: Definition of the Optimization Problem}
The user has to collect all the information about the optimization problem in a single file that has one return argument which is a data structure. The user must define lower and upper variable bounds for \textit{all}  variables and must also define which variables are continuous and integer, respectively.  The objective function must be defined by a function handle and must accept the variable input vector $\mathbf{x}$ and return a scalar objective function value.  The code package and code manual contain detailed examples of how to define each of these elements. 

\subsection{IA-5: Initial Experimental Design}
Each optimization routine (continuous, integer, mixed-integer) starts with an initial experimental design. The user can define the size of the design (IA-6), but the  minimum size depends on the desired surrogate model (IA-3).  Integrality constraints, if present,  are satisfied by definition of the points in  the initial experimental design. The user can choose between three design strategies, namely MATLAB's Latin hypercube design  with 'maximin' option, a symmetric Latin hypercube design~\cite{Ye2000}, and a design that uses  (a subset of) the corner points of the hypercube defined by the variables' upper and lower bounds, and  the center point of the hypercube.  The user may also define starting points to be added to the initial experimental design (IA-7). Next, the expensive objective function is evaluated at the chosen points in the design. This can be done in parallel if the user has the necessary MATLAB Parallel Computing Toolbox installed.\\

\subsection{IA-3: The Surrogate Model}\label{sec:sampling}
The data from the initial experimental design is  used for fitting the surrogate model. The user can choose between various surrogate models and model mixtures (see Table~\ref{tab:settings}).  If a mixture model is used, Dempster-Shafer theory is applied to determine the weights of each model in the mixture~\cite{Muller2010} in each iteration.

\subsection{IA-4: Selecting the Next Sample Point(s)}

The user can define the desired number of points to be selected in each iteration for doing the expensive function evaluation (IA-8), and although any number of points can be used, we recommend  using a multiple of the available processors  in order to best exploit the computing resources if doing synchronous parallel evaluations is an option (MATLAB Parallel Computing Toolbox is required).\\

MATSuMoTo iterates through sample point selection  and updating of the response surface. If  MATSuMoTo is not able to improve the best point found so far in a given number of consecutive trials, the algorithm will start from scratch, i.e.\ MATSuMoTo restarts from creating a new initial experimental design. In this case, the points that have been sampled before  are not taken into account when fitting the surrogate model  in order to avoid getting trapped in the same local minimum. Note that if the problem dimension is large and the number of allowed function evaluations comparatively low, MATSuMoTo may not restart.\\

\subsubsection{Randomized Sample Point Selection}
When using the randomized sample point selection,  the user can choose between a local (\textsf{CANDloc}) and a global (\textsf{CANDglob}) search method. In the local search method a set of candidate points for the next evaluation is created by perturbing the best point found so far. The global search method contains in addition to these candidates a set of points that is  uniformly selected from the whole variable domain. When creating  candidates, the integrality constraints (if there are any) are taken into account, and therefore only candidates that satisfy the integrality constraints are generated.\\

Since the optimization problems are formulated as minimization problems, the point with the lowest function value is the best point found so far, i.e., $\mathbf{x}_{\text{best}} = \text{argmin}_{\mathbf{x}\in\mathcal{D}} f(\mathbf{x})$. The variables of this point are perturbed  with probability
\begin{equation}
P= \left\{
\begin{array}{ll}
	\max\{0.1, 5/d\} & \text{ if  $d > 5$}\\
	1 & \text{else}
\end{array}.
\right.\label{eq:pertprob}
\end{equation}

The perturbation for the local search  is done by randomly adding the value $r\rho\phi$, where $r=\min_{i\in\{1,\ldots,d_1\}}\{x_i^{\text{up}}-x_i^{\text{low}}\}$ and $\phi\sim\mathcal{N}(0,1)$, and where $\rho$ is randomly selected from the set $\{0.2,0.1,0.05\}$ while guaranteeing that the variables' upper and lower bounds are not exceeded.  For \textit{continuous} problems,  each group (points created by perturbation and points created by uniform selection) contains $500d_1$ points. Candidate points that are too close to already evaluated points (closer than 0.1\% of the shortest side of the hyper-rectangle) are discarded to guarantee that the computation of the surrogate model parameters is well-defined. \\

The generation of candidate points for \textit{pure integer} problems is similar as for continuous problems. For the local search, candidates are generated by perturbing the variables of the best point found so far with probability $P$ defined in~(\ref{eq:pertprob}) and $\rho$ is randomly selected from the set $\{1,2,3\}$ . When adding the random perturbations, it is ensured that the integrality constraints are satisfied. Similarly, for the global search, when uniformly selecting candidates from the whole variable domain, the integer constraints must be satisfied, and hence all candidates comply by generation with the integer constraints. For pure integer problems, each group of candidates contains $500d_2$ points. Candidate points that coincide with already evaluated points are discarded. \\

In the local search for \textit{mixed-integer} problems, three groups of candidate points are created by perturbing the best point found so far as follows. The first group contains points that are  generated by perturbing only (a subset of) the continuous variables of $\mathbf{x}_{\text{best}}$ using $\rho$ randomly selected from the set $\{0.2,0.1,0.05\}$. The second group contains points that are  generated by perturbing only (a subset of) the integer variables of $\mathbf{x}_{\text{best}}$ using $\rho$ randomly selected from the set $\{1,2,3\}$. The third group contains points that are  generated by perturbing (a subset of) all variables of $\mathbf{x}_{\text{best}}$. The perturbation probability $P$ of each variable is defined in equation~(\ref{eq:pertprob}). If the global search is desired, a fourth group of candidates is created  by uniformly selecting points from the whole variable domain. All points satisfy the integrality constraints and are within the variables' upper and lower bounds. For each group $125(d_1+d_2)$ points are created and candidate points that are too close to already evaluated points (closer than 0.1\% of the shortest side of the hyper-rectangle) are discarded.\\

In order to select the most promising candidate points for doing the expensive function evaluation, a weighted score  is computed. First, the surrogate model  is used to predict the objective function values at the candidate points. These values are scaled to [0,1] (response surface criterion~\cite{Regis2007b}). Second, the distance of each point to the set of already evaluated points is determined.  Also these values are scaled to [0,1] (distance criterion~\cite{Regis2007b}).   The weight $w_R$  for the response surface criterion cycles through the pattern $\mathcal{C}=<1, 0.75, 0.5, 0.25, 0>$. The weight for the distance criterion  is $w_D=1-w_R$. If $w_R=1$, then the point with the lowest predicted objective function value is selected. In this case the search is local because points that are close to the best point found so far typically have the lowest predicted function values. On the other hand, if $w_R=0$, then the point with the largest distance to the set of already sampled points is chosen. This corresponds to a global search because the chosen point will be in a relatively unexplored area of the variable domain. By cycling through the pattern $\mathcal{C}$, it is possible to repeatedly transition from a local to a global search.  \\

The perturbation range $\rho$ is dynamically adjusted throughout the algorithm. There are counters for tracking how many consecutive  improvement trials failed and succeeded, respectively. If the number of consecutively failed improvement trials exceeds the threshold $\max\{5, d\}$, where $d$ is the problem dimension, the perturbation range is halved. If the number of consecutive successful improvement trials exceeds 3, the perturbation range is doubled. If the number of perturbation range reductions exceeds the value 5, the algorithm restarts from scratch.

\subsubsection{Minimum Point of the Surrogate Model as Sample Point}
MATSuMoTo contains a second sampling strategy that uses the (local) minimum point of the surrogate model  (option \textsf{SurfMin})  in every iteration as new sample point. The surrogate model is in many cases multimodal, i.e.\ there are several local and possibly several global minima present.  Finding the minimum of the surrogate model is compared to the computationally expensive objective function cheap.\\

For continuous optimization problems, MATSuMoTo uses MATLAB's function \textsf{fmincon} for finding the minimum of the surrogate model starting from a randomly selected point in the variable domain. The default settings for  \textsf{fmincon}'s  tuning parameters are used. For integer and mixed-integer problems, MATSuMoTo uses MATLAB's genetic algorithm  \textsf{ga}, which is able to deal with integrality constraints. An initial guess for starting the search is not needed for \textsf{ga} and MATLAB's default genetic algorithm settings are used. \\

In case the surrogate model is unimodal or if the found minimum point is too close to already evaluated points, MATSuMoTo selects the point  that maximizes the minimum distance to all already sampled points  as new evaluation point in order to improve the fit of the surrogate model in rather unexplored regions of the variable domain. It is not necessary to find in every iteration the global minimum of the surrogate model. The surrogate model is an approximation of the true objective function and  minima of the surrogate model do not necessarily have to coincide with the minima of the true objective function. Hence, it is desirable to explore the various local minima of the surrogate model. In order to overcome the possibility of getting trapped in a local optimum, maximizing the minimum distance to already evaluated points helps to improve the global fit of the surrogate model and to  detect new valleys of the objective function landscape where the true global minimum may be located.

\subsubsection{Parallelism}
In each iteration more than one new sample point can be selected and  MATLAB's Parallel Computing Toolbox should be installed when doing so. If the toolbox is not installed, the algorithm performs better when  selecting only a single point in every iteration. For parallel evaluations, a pool of MATLAB workers is started automatically. In order to fully exploit the parallelism, the number of points evaluated in each iteration should not exceed the number of available workers.    Moreover, if the MATLAB Parallel Computing Toolbox is installed, MATSuMoTo performs the function evaluations of the points in the initial experimental design simultaneously.

\subsubsection{Convergence}
For both sampling strategies, it can be shown that MATSuMoTo is asymptotically complete, i.e.\ assuming an indefinitely long run-time and exact calculations, MATSuMoTo will find the global optimum with probability one. This means that MATSuMoTo is able to achieve the goal of not getting trapped in a local optimum and continue the search for improved solutions globally. In practice, however, the number of function evaluations is oftentimes restricted to few hundreds due to the computational expense associated with a function evaluation.

\section{Discussion}\label{sec:discussion}
In this paper we described MATSuMoTo, the MATLAB Surrogate Model Toolbox,  for box-constrained, computationally expensive, black-box, global optimization problems. MATSuMoTo should be used when a single function evaluation is computationally expensive (several minutes or more) because for these problems the computational cost of building the surrogate model is negligible  compared to the computational cost for evaluating the objective function. For problems where the objective function can be computed in a fraction of a second, MATSuMoTo should not be used because for these problems the computational effort for building the surrogate model dominates. \\

MATSuMoTo is to our best knowledge the only surrogate model toolbox that is able to solve continuous, pure integer, and mixed-integer problems. Moreover,  MATSuMoTo allows the user to choose  from a large variety of surrogate models and surrogate model ensembles, three  initial experimental design strategies, and   three sampling strategies. MATSuMoTo can do several function evaluations synchronously, and thus more than one point can be selected for evaluation in each iteration. MATSuMoTo is freely available and is accompanied by a detailed code manual. MATSuMoTo will be extended in the future and will include  features for handling black-box constraints, noisy objective functions, multi-level optimization problems, and multi-objective optimization problems.




\bibliographystyle{plain}
\bibliography{biblio}




%
%

%

\end{document}